\documentclass{amsart}

\usepackage{amssymb}
\usepackage{latexsym}
\usepackage[cmtip,arrow]{xy}
\usepackage{pb-diagram,pb-xy}

\newcommand{\I}{{i}}
\newcommand{\E}{{e}}
\newcommand{\bbZ}{{\mathbb{Z}}}
\newcommand{\bbC}{{\mathbb{C}}}
\newcommand{\bbR}{{\mathbb{R}}}
\newcommand{\bbT}{{\mathbb{T}}}

\newcommand{\Abas}{{a}}
\newcommand{\Dbas}{{d}}
\newcommand{\ebas}{{e}}
\newcommand{\Sbas}{{s}}

\DeclareMathOperator{\Hom}{Hom}
\DeclareMathOperator{\id}{id}
\DeclareMathOperator{\op}{op}
\DeclareMathOperator{\Span}{Span}
\DeclareMathOperator{\lcm}{lcm}

\DeclareMathOperator{\Ann}{Ann}
\newcommand{\bDelta}{\mbox {\boldmath $\Delta$}}

\newcommand{\bS}{\mbox {\boldmath $S$}}
\newcommand{\bepsilon}{\mbox{\boldmath $\epsilon$}}
\newcommand{\Bimod}[1]{{\boldsymbol #1}} 
\newcommand{\Pairing}[2]{{\langle #1 , #2 \rangle}} 
\DeclareMathOperator{\Motimes}{{\otimes_{\Bimod{\Delta}}}} 

\newtheorem{Theorem}{Theorem}
\newtheorem{Proposition}{Proposition}
\newtheorem{Definition}{Definition}
\newtheorem{Corollary}{Corollary}

\newtheorem{Lemma}{Lemma}

\theoremstyle{definition}
\newtheorem{Remark}{Remark}

\begin{document}

\title[Hopfish structure and modules over irrational rotation algebras]{Hopfish structure and modules\\over irrational rotation algebras}

\author{Christian Blohmann}

\address{Department of Mathematics, University of California, Berkeley, CA 94720, USA.\vspace*{-2ex}}
\address{Jacobs University Bremen, Campus Ring 1, 28759 Bremen, Germany.}

\email{blohmann@math.berkeley.edu}

\author{Xiang Tang}

\address{Department of Mathematics, Washington University, St. Louis, MO 63139, USA.}

\email{xtang@math.wustl.edu}

\author{Alan Weinstein}

\address{Department of Mathematics, University of California, Berkeley, CA 94720, USA.}

\email{alanw@math.berkeley.edu}

\thanks{Christian Blohmann was supported by a Marie Curie Fellowship of the European Union; Alan Weinstein was partially supported by NSF grant DMS-0204100.}

\date{Version: April 12, 2006}


\begin{abstract}
Inspired by the group structure on $S^1/ \bbZ$, we introduce a
weak hopfish structure on an irrational rotation algebra $A$ of finite Fourier series. We
consider a class of simple $A$-modules defined by invertible
elements, and we compute the tensor product between these modules
defined by the hopfish structure. This class of
simple modules turns out to generate an interesting commutative unital ring.
\end{abstract}

\maketitle


\section{Introduction}

The starting point for this work is the following two principles.

\begin{enumerate}

\item
If $H$ is a group, many algebras of functions on $H$
(under pointwise multiplication) are Hopf algebras with the coproduct
$\Delta(a)(g,h)=a(gh)$.

\item
If a group $K$ acts on a space $X$, then an algebra which is a crossed
product of $K$ acting on an algebra of functions on $X$ is a good
substitute for an algebra of functions on $X/K$ when $X/K$
is badly behaved as a topological space.

\end{enumerate}
Now let $X=S^1$  be the unit
circle of complex numbers, or $U(1)$, and let $K=\bbZ$ be the subgroup
generated by an element $e^{i\lambda}$, where $\lambda$ is an
irrational multiple of $2\pi$.  Since  the quotient space $H=S^1/\bbZ$ is a group,
the two principles above
suggest that the
crossed product built from $\bbZ$
acting on an algebra of functions on $S^1$ via rotation through the
angle $\lambda$, which we will call irrational rotation algebra, should be something like a Hopf algebra.
But it is well known that an irrational rotation algebra is not a Hopf
algebra; it does not even admit a counit, i.e. a unital
homomorphism to  $\bbC$.  (All of our function spaces will be complex valued.)
As a remedy, the notion of  \textbf{hopfish
algebra} was introduced in \cite{ta-we-zh:hopfish}, based on the
subtext of item 2 above to the effect that the substitution of
a crossed product algebra for an algebra of functions on a
quotient is based on the theory of Morita equivalence, in which
bimodules are interpreted as generalized homomorphisms between
algebras.  Thus, the coproduct and counit of a hopfish algebra $A$
are taken to be $(A\otimes A,A)$- and $(\bbC,A)$-bimodules rather
than algebra homomorphisms.  The antipode is still an ordinary
antihomomorphism from $A$ to itself; we defer the precise
definition to Section 3.

In \cite{ta-we-zh:hopfish} we stopped short of dealing with irrational rotation
algebras, since it appeared that one would need to work in the
world of $C^*$-algebras, which was beyond the purely algebraic
scope of that
paper.  In fact, an algebraic treatment turns out to  be
possible and quite interesting, even if it may not be optimal.  In
the present article, we discuss the hopfish structure of a crossed
product $*$-algebra built from the irrational rotation action of
the integers on the algebra of ``polynomial'' functions on the
circle (i.e. on the algebra of finite Fourier series under pointwise
multiplication).  Actually, the structure is not quite hopfish--we
need a slight weakening of the antipode criterion, but otherwise,
everything works as in \cite{ta-we-zh:hopfish}.

The notion of hopfish algebra is a dualization of that of ``group
object in the category of differentiable stacks''.  In the
geometric language of \cite{be:stack}, \cite{be-xu:stack}, \cite{mr:stability}, and
\cite{tz:stack}, a stack is represented by a groupoid (with
equivalent groupoids representing the same stack), and a morphism
between stacks is represented by a groupoid bibundle. Thus,
the group $S^1/\bbZ$ is considered as a stack represented by the
action groupoid $G=\bbZ \times S^1$. Ignoring that the quotient topology
on the set $S^1/\bbZ$ is trivial, we may construct an equivalence bibundle
between $S^1/\bbZ$ and $G$. Formal composition with this bibundle
and its inverse turns the ordinary group multiplication map from
$S^1/\bbZ \times S^1/\bbZ$ to $S^1/ \bbZ$ into a perfectly
good $(G\times G,G)$-bibundle $B$ representing the product
operation on the stack represented by $G$.  When we dualize, $G$
is replaced by the groupoid algebra $A$, which is just the
irrational rotation algebra. $G\times G$ is then replaced by
$A\otimes A$, and $B$ becomes an $(A\otimes A,A)$-bimodule which
is the coproduct of our hopfish algebra structure on $A$.  (All unsubscripted tensor products are taken over $\bbC$.)

The construction above is presented in detail in Section 2, along
with a construction for the counit of $A$, derived in a similar
way from the inclusion of the identity element into $S^1/\bbZ.$ In
Section 3.2, we discuss the antipode of $A$, which should be
derived from the inversion map on $S^1/\bbZ$.  Here it turns out
that the  object which we construct does not quite satisfy the definition
in \cite{ta-we-zh:hopfish} because of difficulties with the
duality of infinite-dimensional vector spaces. Therefore, we
weaken the notion of antipode to accommodate this example.

In the last part of the paper, we use the coproduct bimodule to
construct a tensor product operation on the collection of isomorphism
classes of $A$-modules, and we study the behavior of this application
when applied to a nice class of cyclic modules.  From the resulting
algebraic structure, it is possible to reconstruct the quotient group
$S^1/\bbZ$.

It is important to note that, although the irrational rotation
algebra may be viewed as a deformation of the algebra of functions
on a 2-dimensional torus \cite{ri:quantization}, our hopfish
structure is {\em not} a deformation of the Hopf structure
associated with the group structure on the torus.  Rather, the
classical limit of our hopfish structure is a second symplectic
groupoid structure on $T^*\bbT^2$ which is compatible with the one
described in \cite{we:rotation}, whose quantization is the
multiplication in the irrational rotation algebra.   We thus seem
to have a symplectic double groupoid which does not arise from a
Poisson Lie group as do those in \cite{lu-we:groupoides} and
\cite{za:quantum}.  It is possible that such double groupoids are
in general the classical limit of hopfish algebras and thus
represent a useful generalization of Poisson Lie groups. We would
like to also mention that the hopfish algebra structure on an
irrational rotation algebra is closely related to the para-Hopf
algebroid structure introduced in \cite{ka-ra:para-hopf}. In
general, there seem to be interesting connections between the
notions of hopfish algebra and Hopf algebroid. We hope to pursue
these issues in the future.

From the viewpoint of higher algebra, a groupoid with compatible
group structure can be viewed as a special case of a 2-groupoid
over a point \cite{ba-la:2-groups}. In this language the groupoid
multiplication on $S^1 \times \bbZ$ is the vertical composition of
2-morphisms, whereas the product group multiplication on $S^1
\times \bbZ$ is the additional horizontal composition.
Accordingly, on the space of functions on $S^1 \times \bbZ$ both
compositions yield convolution products and both inverses yield
star structures. The vertical convolution product is the usual
convolution of the groupoid algebra which encodes the ``bad''
topology of the stack it represents. The horizontal convolution is
the group convolution which encodes the additional group
structure. A preliminary study of this approach for the irrational
rotation algebra is very promising: The group convolution
operation on the states corresponding to the cyclic modules we
consider here is closely related to the tensor product obtained by
the hopfish structure. Moreover, the weak hopfish antipode turns
out to be the composition of the vertical and the horizontal star
involution. This is intriguing and suggests further research.

There are several other important steps in the study of hopfish algebras
which we have not yet taken.  We still work in the category in
which objects are algebras and morphisms are isomorphism classes
of bimodules, rather than in the 2-category in which morphisms are
bimodules and 2-morphisms are bimodule isomorphisms.  As a result,
we do not make the modules over $A$ into a tensor category, but
simply work with the isomorphism classes of these modules.
Second, we are still working purely algebraically, rather than in the
world of $C^*$-algebras, where the use of topological tensor
products may allow us to circumvent the weak property of the
antipode.  Here we should note, though, that even for Hopf
algebras, where the coproduct is a homomorphism rather than a
bimodule, the extension of the theory to include topologies is
highly nontrivial, e.g. \cite{bo-fl-ge-pi:hidden,wo:quant-gp}.

\bigskip
\noindent {\bf Acknowledgements.}  We would like to thank Marc
Rieffel, Noam Shomron, and Chenchang Zhu for helpful discussions, and Yvette Kosmann-Schwarzbach for comments on the manuscript.

\section{Construction of the hopfish algebra}
\label{sec-construction}
We begin by constructing for the action groupoid $S^1\times \bbZ$
the bibundles which encode the
group multiplication, unit, and inversion of the quotient group
$S^1/\bbZ$.  Spaces of functions on these bibundles are the
bimodules giving the (weak) hopfish structure on the irrational
rotation algebra.

\subsection{The action groupoid and its groupoid algebra}
\label{subsec-action}
As a set, the action groupoid is\footnote{We
 follow the convention of operator algebra theory, in which the group
 in a crossed product is written to the right, even though it is
 acting from the left.} $G = S^1 \times \bbZ$.  The left and
right projections onto the base $S^1$ are defined by
$\alpha(\theta,k) = \theta + \lambda k$ and
$\beta(\theta,k) = \theta$.   If we think of a groupoid element
as an arrow pointing from an
element $\theta$ of $S^1$ to its image under the action of $k \in
\bbZ$, then $\beta$ is the source map and $\alpha$ the target.
When
$\beta(\theta_1,k_1) = \alpha(\theta_2,k_2)$, that is, when
 $\theta_1 = \theta_2 +
\lambda k_2$, the product of the two groupoid elements is
$(\theta_2, k_1 + k_2 )$. The identity bisection is given by the
natural embedding $S^1 \hookrightarrow S^1 \times \bbZ$, $\theta
\mapsto (\theta, 0)$.  Inversion is the mapping
$(\theta,k)\mapsto (\theta -\lambda k, -k)$.

With the product topology and
differentiable structure, $G=S^1 \times \bbZ$
is a smooth, \'etale
groupoid. The family of counting measures on the source fibres
is a natural Haar system,
leading to an associative convolution product on the space
$C^\infty_\mathrm{c}(G)$ of smooth, compactly supported
functions, defined by

\begin{equation}
\label{eq:Convolution2}
 (a * b)(\theta,k)
 := \sum_{k' \in \bbZ} a(\theta + \lambda k',k-k') b(\theta,k') \,.
\end{equation}
For the functions
\[
 \Abas_{nl} \in C(S^1 \times \mathbb{Z})
 \,,\quad \Abas_{nl}(\theta, k)
 := \E^{\I n \theta} \delta_{lk}
 \,,\quad n,l \in \mathbb{Z} \,,
\]
the convolution product~\eqref{eq:Convolution2}  is
given by the simple formula
\begin{equation}
\label{eq:AProduct}
 \Abas_{n_1 l_1} * \Abas_{n_2 l_2}
 = \E^{\I \lambda n_1 l_2} \Abas_{n_1 + n_2, l_1 + l_2} \,.
\end{equation}
Furthermore, we have a conjugate linear involution  defined by $a^*(g) = \overline{a(g^{-1})},$ which yields
\begin{equation}
\label{eq:AStar}
 \Abas_{nl}^* = \E^{\I \lambda n l } \Abas_{-n,-l} \,.
\end{equation}
Note, that all $\Abas_{nl}$ are unitary and that $\Abas_{00}$ is the
unit element. The vector space spanned by $\{ \Abas_{nl} \,|\, n,l
\in \bbZ \}$ with product~\eqref{eq:AProduct} and
involution~\eqref{eq:AStar} is the unital $*$-algebra generated by
$\Abas_{10}$ and $\Abas_{01}$.
It is on this algebra, dense in the irrational
rotation $C^*$-algebra, that we will focus attention in this paper.
Denoting it by $A$, we may think of it as the algebra of ``regular
functions on an algebraic quantum torus.''

\subsection{The coproduct}
We wish to construct bibundles which represent for the groupoid
$S^1\times \bbZ$ the mappings representing the group structure on $S^1/\bbZ$.
These are obtained by composing the mappings with the bibundle (and
its inverse) representing the equivalence between the action groupoid
and the quotient space.  Even though the quotient space has a ``bad
topology,'' the final bibundles will be perfectly nice.

The general scheme for our constructions is the following.
Let $G$ be a groupoid over a set $X$, $H$ a groupoid over $Y$, and $f: X/G
\rightarrow Y/H$ a morphism between the quotients.  Composing the
bibundle representing $f$ with the bibundles representing the
equivalences between the groupoids and their quotient spaces amounts
to filling out the following array of pullback diagrams, working from
the given diagrams at the bottom to the final diagram
$X \leftarrow X \times_{Y/H} Y \rightarrow Y.$

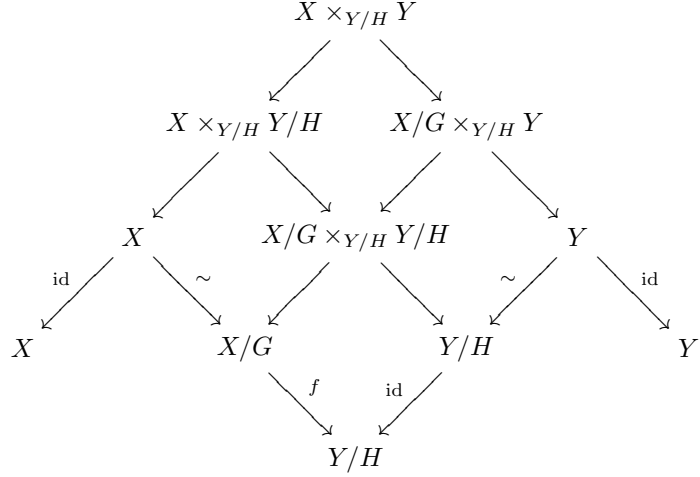
\[
\dgARROWLENGTH=0.2em
\dgHORIZPAD=-2.5em
\begin{diagram}
\node[4]{X \times_{Y/H} Y}
 \arrow{ssw,l}{}
 \arrow{sse,t}{}
\\[2]
\node[3]{X  \times_{Y/H} Y/H}
 \arrow{ssw,l}{}
 \arrow{sse,t}{}
\node[2]{X/G \times_{Y/H} Y}
 \arrow{ssw,l}{}
 \arrow{sse,t}{}
\\[2]
\node[2]{\qquad X \qquad}
 \arrow{ssw,l}{\id}
 \arrow{sse,t}{\sim}
\node[2]{X/G \times_{Y/H} Y/H}
 \arrow{ssw,l}{}
 \arrow{sse,t}{}
\node[2]{\qquad Y \qquad}
 \arrow{ssw,l}{\sim}
 \arrow{sse,t}{\id}
\\[2]
\node{\qquad X \qquad}
\node[2]{\qquad X/G \qquad}
 \arrow{sse,t}{f}
\node[2]{\qquad Y/H \qquad}
 \arrow{ssw,l}{\id}
\node[2]{\qquad Y \qquad}
\\[2]
\node[4]{\qquad Y/H \qquad}
\end{diagram}
\]

Now let $\pi : S^1 \rightarrow S^1/\bbZ$ denote the canonical epimorphism
and $+$ the group operation on $S^1/\bbZ$. We have the following
commutative diagram of two nested pull-back squares, which
is a subdiagram of the previous one, tilted by $45$ degrees.

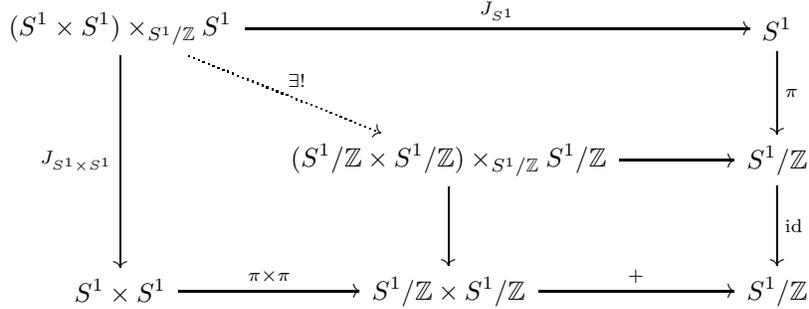
\[
\begin{diagram}
\node{(S^1 \times S^1)\times_{S^1/\bbZ} S^1}
 \arrow[2]{s,l}{J_{S^1\times S^1}}
 \arrow[2]{e,t}{J_{S^1}}
 \arrow{se,t,..}{\exists!}
\node[2]{S^1}
 \arrow{s,r}{\pi}
\\
\node[2]{(S^1/\bbZ \times S^1/\bbZ) \times_{S^1/\bbZ} S^1/\bbZ}
 \arrow{s}
 \arrow{e}
\node{S^1/\bbZ}
 \arrow{s,r}{\id}
\\
\node{S^1 \times S^1}
 \arrow{e,t}{\pi \times \pi}
\node{S^1/\bbZ \times S^1/\bbZ}
 \arrow{e,t}{+}
\node{S^1/\bbZ}
\end{diagram}
\]
The inner pull-back $(S^1/\bbZ \times S^1/\bbZ) \times_{S^1/\bbZ} S^1/\bbZ$ is simply the graph of addition on $S^1/\bbZ$, the pull-back projections being the range and image maps. Any object in the left upper corner which makes the diagram commute can be viewed as a lift of $\mathrm{graph}(+)$ to $S^1$. The diagonal arrow is the unique map which exists by the universal property of the inner pull-back. For example, the graph of addition on $S^1$ is such a lift of $\mathrm{graph}(+)$. The left upper corner of the outer pull-back square can then be viewed as the universal lift into which all other lifts map uniquely.

Let us determine this universal lift explicitly. Denote the image of an element $\theta \in S^1$ under the canonical epimorphism by $[\theta]$. The pull-back is the set
\[
 (S^1 \times S^1)\times_{S^1/\bbZ} S^1 =
 \{(\theta_1,\theta_2, \theta)
 \in (S^1 \times S^1) \times S^1
 \,|\, [\theta_1] + [\theta_2] = [\theta] \}
\]
together with the pull-back projections
\[
 J_{S^1 \times S^1}(\theta_1,\theta_2, \theta) = (\theta_1,\theta_2)
 \,,\qquad
 J_{S^1}(\theta_1,\theta_2, \theta) = \theta \,.
\]
It is convenient to identify the pull-back as a set with $S^1 \times S^1 \times \bbZ$ by the map
\[
\begin{aligned}
 (S^1 \times S^1)\times_{S^1/\bbZ} S^1
 &\stackrel{\cong}{\longrightarrow} S^1 \times S^1 \times \mathbb{Z}\\
 (\theta_1, \theta_2, \theta)
 &\longmapsto (\theta_1, \theta_2,(\theta - \theta_1 - \theta_2)/\lambda) \,.
\end{aligned}
\]
For an element $(\theta_1,\theta_2,k) \in S^1 \times S^1 \times \bbZ$ the pull-back projections are
\[
 J_{S^1 \times S^1}(\theta_1,\theta_2, k) = (\theta_1,\theta_2)
 \,,\qquad
 J_{S^1}(\theta_1,\theta_2, k) = \theta \,.
\]
From the left and right groupoid actions of $G$ on $S^1$ the
pull-back inherits a left action of $G\times G$ and a right action
of $G$, the pull-back projections being the momentum maps of these
actions. Explicitly, the left action of $(\phi_1,l_1,\phi_2,l_2)
\in G \times G$ reads
\[
(\phi_1,l_1,\phi_2,l_2) \cdot (\theta_1, \theta_2, k) = (\theta_1
+ \lambda l_1, \theta_2 + \lambda l_2, k - l_1 - l_2),
\]
which is defined if $\phi_1 = \theta_1$ and $\phi_2 = \theta_2$. The right action of $(\phi,l) \in G$ reads
\[
 (\theta_1, \theta_2, k) \cdot (\phi, l)
 = (\theta_1, \theta_2, k - l)
 \,,
\]
which is defined if $\theta_1 + \theta_2 + \lambda k = \phi +
\lambda l$. Together with these actions the pull-back becomes a
groupoid bibundle.

The space of smooth functions on this groupoid bibundle can now be
equipped with the structure of a bimodule over the groupoid
algebras: Let $m \in C^{\infty}(S^1 \times S^1 \times \bbZ)$, let
$a \in A \subset C^\infty_\mathrm{c}(G)$, and let $x \in S^1
\times S^1 \times \bbZ$. Then
\begin{equation}
\label{eq:rightact1}
 (m \cdot a)(x)
 := \sum_{g \in \beta^{-1}\left(J_{S^1}(x)\right)}
 m( x \cdot g^{-1} ) \,
 a(g)
\end{equation}
defines a right action of the groupoid algebra $A$.

The analogous construction on the left side yields a left $A
\otimes A$-action such that the space of functions on the bibundle
becomes an $(A\otimes A,A)$-bimodule. Again, we are interested
only in the algebraic picture so we choose a set of functions
\begin{equation}
\label{eq:Deltabasis}
 \Dbas_{n_1 n_2 l} \in C(S^1 \times S^1 \times \bbZ)
 \,,\quad \Dbas_{n_1 n_2 l}(\theta_1, \theta_2, k)
 := \E^{\I n_1\theta_1} \E^{\I n_2\theta_2} \delta_{lk}
 \,,\quad n_1,n_2, l \in \bbZ \,.
\end{equation}
The right action~\eqref{eq:rightact1} of the basis of $A$ on these
function is easily computed to be
\begin{equation}
\label{eq:AactDelta}
 \Dbas_{n_1 n_2 l} \cdot \Abas_{mj}
 = \E^{\I \lambda m (l-j)} \Dbas_{n_1+m, n_2+m, l-j} \,.
\end{equation}
For the left action of $A \otimes A$ we obtain
\begin{equation}
\label{eq:bim-coprod}
 (\Abas_{m_1 j_1} \otimes \Abas_{m_2 j_2})
 \cdot \Dbas_{n_1 n_2 l}
 = \E^{-\I \lambda [ (n_1 + m_1) j_1 + (n_2 + m_2) j_2]}
 \Dbas_{n_1+m_1, n_2+m_2, l - j_1 - j_2} \,.
\end{equation}
The $(A\otimes A, A)$-bimodule  $\bDelta$ spanned by
$\{ \Dbas_{n_1 n_2 l} \}$ is a natural candidate
for the hopfish coproduct on $A$.

\subsection{The counit}

The unit element of the group $S^1/\bbZ$ can be viewed as a map $e
: \{ \mathrm{pt} \} \rightarrow S^1/\bbZ$, $e(\mathrm{pt}) = 0$.
The pull-back corresponding to the universal lift to $S^1$, is the
set
\[
 \{ (\mathrm{pt}, \theta) \in \{\mathrm{pt}\} \times S^1 \,|\, e(\mathrm{pt}) = [\theta] \}
 \cong \bbZ
\]
with the right pull-back projection $J_{S^1}(k) = \lambda k$. The right action of $(\phi,l) \in G$ is
\[
 k \cdot (\phi,l) = k - l \,,
\]
which is defined if $\lambda k = \phi + \lambda l$.

Once more, we choose a set of functions on this right groupoid bundle,
\begin{equation}
\label{eq:eBasis}
 \ebas_l \in C(\bbZ)
 \,,\quad \ebas_l(k)
 := \delta_{lk}
 \,,\quad l \in \mathbb{Z} \,.
\end{equation}
The right action of $A$ is computed to be
\[
 \ebas_l \cdot \Abas_{mj} = \E^{\I\lambda m (l-j)} \ebas_{l-j} \,.
\]
The right $A$-module $\Bimod{\epsilon}$ spanned by $\{\ebas_l\}$ by
 is the natural candidate for the
hopfish counit of $A$.

\subsection{The antipode}

The pull-back corresponding to the universal lift of the inversion
map on $S^1 / \bbZ$ to $S^1$, is
\[
 \{(\theta_1,\theta_2) \in S^1 \times S^1 \,|\, [-\theta_1] = [\theta_2] \}
 \cong S^1 \times \bbZ \,,
\]
where the identification is $(\theta_1, \theta_2) \mapsto
(-\theta_2, (\theta_1 + \theta_2)/\lambda)$. The pull-back
projections are $J_{S^1}^{\mathrm{left}}(\theta,k) = \theta +
\lambda k$ on the left and $J_{S^1}^{\mathrm{right}}(\theta) =
-\theta$ on the right.  In light of the axioms of a hopfish
algebra, however, we have to view the pull-back as a left $S^1
\times S^1$ bundle with the bundle projection
$J_{S^1}^{\mathrm{left}} \times J_{S^1}^{\mathrm{right}}$. The
corresponding left action of $(\phi_1,l_1,\phi_2,l_2) \in G \times
G$ is
\[
 (\phi_1,l_1,\phi_2,l_2) \cdot (\theta, k)
 = (\theta - \lambda l_2, k + l_1 + l_2) \,,
\]
which is defined if $\phi_1 = \theta + \lambda k$ and $\phi_2 = - \theta$.

As a set, the pull-back of the graph of the inverse is isomorphic
to the groupoid. This suggests choosing as basis for the bimodule
the same set of functions as for the groupoid algebra:
\[
 \Sbas_{nl} \in C(S^1 \times \mathbb{Z})
 \,,\quad \Sbas_{nl}(\theta, k)
 := \E^{\I n \theta} \delta_{lk}
 \,,\quad n,l \in \mathbb{Z} \,.
\]
The left action of $A \otimes A$ on these
functions is computed to be
\begin{equation}
\label{eq:Smodule}
 (\Abas_{m_1 j_1} \otimes \Abas_{m_2 j_2})
 \cdot \Sbas_{n l}
 = \E^{\I \lambda(m_1 l)} \E^{\I(n+m_1-m_2)j_2}
 \Sbas_{n+m_1-m_2, l + j_1 + j_2} \,.
\end{equation}
The left $A \otimes A$-module  $\Bimod{S}$ spanned by $\{\Sbas_{nl}\}$ by
is the natural candidate for the hopfish antipode.

\section{Verification of the axioms}
In this section, we study relations among the bimodules  $(\bDelta, \bepsilon,
\bS)$ defined in Section \ref{sec-construction}.

\subsection{The sesquiunital sesquialgebra}
When algebra homomorphisms are replaced by bimodules, the notion of
unital bialgebra becomes that of sesquiunital sesquialgebra
\cite{ta-we-zh:hopfish}.
\begin{Definition}
\label{dfn-sesqui} A {\bf sesquiunital sesquialgebra} over a
commutative ring $k$ is a unital $k$-algebra $A$ equipped with an
$(A\otimes A, A)$-bimodule $\bDelta$ (the {\bf
 coproduct}) and a $(k,A)$-module (i.e.
 a right $A$-module) $\bepsilon$
(the {\bf counit}), satisfying the following properties.

\begin{enumerate}
\item[\rm{(H1)}](coassociativity)  The
 $(A\otimes A\otimes A,A)$-bimodules
$(A\otimes \bDelta)\otimes_{A\otimes  A}\bDelta$ and
$(\bDelta\otimes A)\otimes_{A\otimes  A}\bDelta$ are isomorphic.
\item[\rm{(H2)}] (counit)  The $(k\otimes A,A)=
(A\otimes  k,A)=(A,A)$-bimodules \\
 $(\bepsilon\otimes  A)\otimes_{A\otimes  A}\bDelta$ and
$(A\otimes_ k \bepsilon)\otimes_{A\otimes A}\bDelta$ are both
 isomorphic to $A$.
\end{enumerate}
\end{Definition}

\begin{Proposition}
Let $A$ be the algebra defined in
Section \ref{subsec-action}.
 The coproduct $\bDelta$ and counit $\bepsilon$ spanned by the
bases in Eq. (\ref{eq:Deltabasis}) and Eq. (\ref{eq:eBasis}) define a
sesquiunital sesquilinear algebra structure on $A$.

\begin{proof} We verify the coassociativity (H1) for $\bDelta$;
the proof for the counit property (H2) is similar.

Since $A$ is free of rank one over itself, the
bimodule $(A\otimes \bDelta)\otimes_{A\otimes A}\bDelta$ is the
linear span of elements of the form
\[
(1\otimes d_{n_1,n_2,l}) \otimes_{\bbC \otimes A} d_{n'_1,n'_2,l'},\ \ \ n_1,
n_2, l, n'_1, n'_2, l'\in \bbZ.
\]
By Eqs. (\ref{eq:AactDelta}) and (\ref{eq:bim-coprod}),
\begin{eqnarray*}
&&e^{i\lambda m(l-j)}(1\otimes d_{n_1+m,n_2+m,l-j})\otimes
d_{n'_1,n'_2,l'} \\
&=&e^{-i\lambda(n'_2+m)j}(1\otimes d_{n_1, n_2, l})\otimes
d_{n'_1, n'_2+m, l'-j}.
\end{eqnarray*}
This relation tells us that the bimodule $(A\otimes
\bDelta)\otimes_{A\otimes A}\bDelta$
is spanned (over $\bbC$) by the elements
\[
(1\otimes d_{n_1, n_2, 0})\otimes d_{n_3,0, l},\ \ \ n_1, n_2,
n_3, l\in \bbZ.
\]
It is easy to see that these generators are linearly independent and
form a basis of $(A\otimes \bDelta)\otimes_{A\otimes A}\bDelta$.

The left $A$-module structure on $(A\otimes
\bDelta)\otimes_{A\otimes A}\bDelta$ is computed as follows.
\begin{eqnarray*}
&&(1\otimes d_{n_1, n_2, 0})\otimes d_{n_3,0, l}\cdot a_{mj}\\
&=&(1\otimes d_{n_1, n_2, 0})\otimes e^{i\lambda m(l-j)}d_{n_3+m,
m, l-j}\\
&=&e^{i\lambda m(l-j)}(1\otimes d_{n_1+m, n_2+m, 0})\otimes
d_{n_3+m, 0, l-j}.
\end{eqnarray*}
And the right $A\otimes A\otimes A$-module structure on $(A\otimes
\bDelta)\otimes_{A\otimes A}\bDelta$ is computed as follows.
\begin{eqnarray*}
&&(a_{m_1,j_1}\otimes a_{m_2, j_2}\otimes a_{m_3,
j_3})\cdot(1\otimes d_{n_1, n_2, 0})\otimes d_{n_3,0, l}\\
&=&e^{-i\lambda\Theta}(1\otimes
d_{n_1+m_2, n_2+m_3, -j_2-j_3})\otimes d_{n_3+m_1,0,l-j_1}\\
&=&e^{-i\lambda\Theta}(1\otimes d_{n_1+m_2, n_2+m_3, 0})\otimes
d_{n_3+m_1, 0,l-j_1-j_2-j_3},
\end{eqnarray*}
where $\Theta=(n_3+m_1)j_1+(n_1+m_2)j_2+(n_2+m_3)j_3$.

A similar computation shows that $(\bDelta\otimes
A)\otimes_{A\otimes A}\bDelta$ has a basis
\[
(d_{n_1, n_2,0}\otimes 1)\otimes d_{0,n_3, l},\ \ \ n_1, n_2, n_3,
l\in \bbZ.
\]
We define the following map $I:(A\otimes \bDelta)\otimes_{A\otimes
A}\bDelta\to (\bDelta\otimes A)\otimes_{A\otimes A}\bDelta$ by
\[
I((1\otimes d_{n_1, n_2, 0})\otimes d_{n_3,0,l})=(d_{n_3,
n_1,0}\otimes 1)\otimes d_{0,n_2,l},\ \ \ n_1, n_2, n_3, l\in
\bbZ.
\]
It is easy to check that $I$ is a bimodule
isomorphism.
\end{proof}
\end{Proposition}

\subsection{The antipode}
We recall from
 \cite{ta-we-zh:hopfish} the definition of an antipode
for a hopfish algebra.

\begin{Definition}\label{dfn-preantipode}
A {\bf preantipode} for a sesquiunital sesquialgebra $A$ over $k$
is a left $A\otimes  A$-module $\bS$ together with an isomorphism
of its $k$-dual with the right $A\otimes A$-module
$\Hom_A(\bepsilon,\bDelta)$ of left $A$-linear maps to $\Bimod{\epsilon}$ from $\Bimod{\Delta}$.\footnote{We use the convention that homomorphisms map from right to left in order to be consistent with \cite{ta-we-zh:hopfish}.}

If a preantipode $\bS$, considered as an $(A, A^{op})$-bimodule,
is a free left $A$-module of rank 1, we call $\bS$ an {\bf
antipode} and say that $A$ along with $\bS$ is a {\bf hopfish
algebra}.
\end{Definition}

The definition of (pre)antipode can reformulated as the
following two conditions.
\begin{enumerate}

\item[(H3)] (preantipode) The dual module $\Bimod{S}^*$ and the space
of right $A$-linear maps $\Hom_A(\Bimod{\epsilon}, \Bimod{\Delta})$
are isomorphic as right $A \otimes A$-modules.

\item[(H4)] (antipode) As a left $A = A \otimes \bbC$-module, $\Bimod{S}$ is
free of rank one.

\end{enumerate}

When the left $A \otimes A$-module $\Bimod{S}$ is viewed as an
$(A,A^{\op})$-bimodule, Axiom (H4)  states that
$\Bimod{S}$ is the modulation of a homomorphism of algebras $S :
A^{\op} \rightarrow A$.  In fact, this is easily verified. From
Eq.~\eqref{eq:Smodule} we get
\[
 (\Abas_{m j} \otimes 1) \cdot \Sbas_{0 0} = \Sbas_{mj} \,,
\]
which shows that $\Sbas_{00}$ is a
basis of $\Bimod{S}$ as an $A\otimes \bbC$-module.
Furthermore,
\[
(1 \otimes \Abas_{m j}) \cdot \Sbas_{0 0}
 = \E^{-\I \lambda m j} \Sbas_{-m,j}
 = (\E^{-\I \lambda m j} \Abas_{-m,j} \otimes 1) \cdot \Sbas_{00} \,,
\]
from which we conclude that $\Bimod{S}$ is (isomorphic to) the
modulation of the homomorphism
\begin{equation}
\label{eq:Smap}
 S : A^{\op} \rightarrow A \,,\qquad
 S(\Abas_{mj}) = \E^{-\I \lambda m j} \Abas_{-m,j} \,.
\end{equation}
Note, that $S^2 = \id$ as would be the case for a cocommutative
Hopf algebra, and that $S \circ * \circ S \circ * = \id$ as
expected. However, it turns out that axiom (H3) does not hold for
$\Bimod{S}$. In fact, we have the following proposition:

\begin{Proposition}
 The sesquiunital sesquialgebra $(A, \Bimod{\epsilon}, \Bimod{\Delta})$ does not admit a hopfish antipode.
\end{Proposition}

\begin{proof}
Any isomorphism $\psi: \Bimod{S}^* \rightarrow
\Hom_A(\Bimod{\epsilon}, \Bimod{\Delta})$ of right $A \otimes
A$-modules is a fortiori an isomorphism of right $A \otimes
1$-modules. Such an isomorphism must map eigenvectors of algebra
elements to eigenvectors.   Since any hopfish antipode $\Bimod{S}$ is
isomorphic to $A$ as a left $A \otimes \bbC = A$-module,  its dual
$\Bimod{S}^*$ is isomorphic to $A^*$ as  a right $A$-module. Consider $z
 \in A^*$ defined by $z(\Abas_{nl}) = \delta_{n,0}$. We have
\begin{equation}
\label{eq:zEigen}
 (z \cdot \Abas_{01} )(\Abas_{nl}) = z(\Abas_{01} * \Abas_{nl})
 = z(\Abas_{n,l+1}) = \delta_{n,0} = z(\Abas_{nl}) \,,
\end{equation}
that is, $z$ is an eigenvector of $\Abas_{01}$. We will now show that this eigenvector cannot be mapped to an eigenvector and conclude that $\psi$ cannot be an isomorphism.

Let us first determine the $A \otimes A$-module structure of $\Hom_A(\Bimod{\epsilon}, \Bimod{\Delta})$ explicitly. Relabel the basis~\eqref{eq:Deltabasis} of $\Bimod{\Delta}$ by
\[
 \tilde{\Dbas}_{n n_2 l} := \Dbas_{n_2 + n, n_2, l}
 \,,\quad n,n_2,l \in \bbZ \,,
\]
that is, we substitute $n := n_1 -  n_2$. The right
action~\eqref{eq:AactDelta} of $A$ now reads
\[
 \tilde{\Dbas}_{n n_2 l} \cdot \Abas_{mj}
 = \E^{\I\lambda m (l-j)} \tilde{\Dbas}_{n, n_2+m,l-j} \,.
\]
 From this we can see that, as a right $A$-module,
 $\Bimod{\Delta}$ is the direct sum of the modules
\[
 \Bimod{\Delta} \cong \bigoplus_{n \in \bbZ} V_{n}
 \,,\quad V_{n} := \Span_\bbC
 \{ \tilde{\Dbas}_{n n_2 l} \,|\, n_2,l \in \bbZ \} \,.
\]
Each $A$-module $V_{n}$ is simple, cyclic, generated by $\tilde{\Dbas}_{n00} = \Dbas_{n00}$, and free. Hence, $\Bimod{\Delta}$ is a free right $A$-module with $A$-basis $D := \{\Dbas_{n00} \,|\, n \in \bbZ \}$. We deduce that we have isomorphisms of vector spaces\footnote{Recall our convention that $\Hom(X,Y)$ denotes the set of homomorphisms to $X$ from $Y$.}
\[
 \Hom_A(\Bimod{\epsilon}, \Bimod{\Delta})
 \cong
 \Hom_{\mathrm{Set}}(\Bimod{\epsilon}, D) \,.
\]
That is, every homomorphism in $\zeta \in \Hom_A(\Bimod{\epsilon},
\Bimod{\Delta})$ is determined by its values on $\Dbas_{n00}$, which
can be chosen freely.  Such a homomorphism
can, therefore, be represented by the matrix elements $\zeta^l_n
\in \bbC$ defined as
\begin{equation}
\label{eq:MatrixEl}
 \zeta(\Dbas_{n00}) = \sum_l \zeta^l_n \ebas_l \,,
\end{equation}
where $\{ \ebas_l \}$ is the basis of $\Bimod{\epsilon}$ defined in
Eq.~\eqref{eq:eBasis}.  The sum over $l$
must be finite for each $n$.

The right action of $a \otimes b \in A \otimes A$ on $\zeta \in \Hom_A(\Bimod{\epsilon}, \Bimod{\Delta})$ is by pullback:
\[
 \bigl( \zeta \cdot (a \otimes b) \bigr)(\Dbas)
 := \zeta\bigl( (a \otimes b) \cdot  \Dbas \bigr)
 \,,\qquad \Dbas \in \Bimod{\Delta} \,.
\]
For the action of the basis of $A \otimes \bbC$, we evaluate
\[
 \bigl( \zeta \cdot (\Abas_{m j} \otimes 1) \bigr) (\Dbas_{n00})
 = \E^{-\I \lambda j(n + m)} \zeta( \Dbas_{n+m, 0, 0}) \cdot \Abas_{0j} \,,
\]
where the action of the right hand side is the right action of $A$ on $\Bimod{\epsilon}$. This can be expressed in terms of the matrix elements as
\begin{equation}
\label{eq:HomAAact}
 \bigl( \zeta \cdot (\Abas_{m j} \otimes 1) \bigr)^l_n
 = \E^{-\I \lambda j(n + m)} \zeta^{l+j}_{n+m} \,.
\end{equation}

Now we come back to the eigenvector $z$ defined in
Eq.~\eqref{eq:zEigen}. Let $\zeta := \psi(z) \in
\Hom_A(\Bimod{\epsilon}, \Bimod{\Delta})$. If $\psi$ is right
$A\otimes 1$-linear we have $\zeta \cdot (\Abas_{01} \otimes 1) =
\zeta$. Assume that there is a non-zero matrix element
$\zeta^l_n$.  By Eq.~\eqref{eq:HomAAact}, we have $\zeta^l_n =
(\zeta \cdot (e_{01} \otimes 1))^l_n = \exp(-\lambda \I n)
\zeta^{l+1}_n$. It follows by induction that $\zeta^{k}_n$ is nonzero
for all $k \geq l$. But this contradicts the fact that for fixed $n$ the
sum $\sum \zeta^l_n \ebas_l \in \Bimod{\epsilon}$ is finite. Hence,
$\zeta$ must be zero, so the kernel of $\psi$ is not empty.
\end{proof}

\subsection{The weak hopfish antipode}

The nonexistence of a hopfish antipode is tied to the algebraic dual
appearing in the axiom (H3) for the preantipode.  For infinite
dimensional vector spaces, forming the algebraic dual is an inconvenient
operation, as it raises the dimension (i.e. the cardinality of a
basis).   Thus, infinite dimensional vector spaces are
never reflexive, and countably infinite spaces do not have a predual at
all. This suggests replacing axiom (H3) by a weaker notion of
preantipode which does not involve the dual.

Such a weaker definition of preantipode is obtained by substituting
the notion of duality in the definition of the preantipode with that
of a dual pairing. Recall that a dual pairing of complex vector
spaces $U$ and $V$ is a $\bbC$-bilinear map $\Pairing{~}{~}: U \times
V \rightarrow \mathbb{C}$. Let $A$ be a unital ring. If $U$ is a right
$A$-module and $V$ is a left $A$-module, the pairing is called
$A$-tensorial if $\Pairing{u \cdot a}{v} = \Pairing{u}{a \cdot v}$ for
all $u \in U$, $v \in V$, and $a \in A$.  If, for each $u\in U$,
 the
vanishing of $\Pairing{u}{v} = 0$ for all $v\in V$ implies $u =0$, then the pairing is
called non-degenerate in $U$.  Nondegeneracy in $V$ is defined
similarly, and the pairing is simply called non-degenerate if
it is non-degenerate in both $U$ and $V$.


\begin{Definition}
Let $\Bimod{S}$ be a left $A \otimes A$-module. $\Bimod{S}$ is called a weak hopfish preantipode, if the following holds:
\begin{enumerate}
\item[{\rm (H3')}] There is a non-degenerate $A \otimes A$-tensorial dual pairing of $\Hom_A(\Bimod{\epsilon}, \Bimod{\Delta})$ and $\Bimod{S}$.
\end{enumerate}
A weak hopfish preantipode is called a weak hopfish antipode if it satisfies axiom \textup{(H4)}.
\end{Definition}

This is really a weaker notion of antipode, since, with respect to the
canonical pairing of a module with its algebraic dual, every hopfish
preantipode is a weak preantipode and every hopfish antipode a weak
hopfish antipode.

In the case of the irrational rotation algebra, we have already seen
that $\Bimod{S}$ satisfies axiom (H4). It is isomorphic to the
modulation of the homomorphism $S$ defined in Eq.~\eqref{eq:Smap}, the
isomorphism being induced by identifying $\Sbas_{00}$ with $\Abas_{00}
= 1_A$. Hence, for every element $s \in \Bimod{S}$ there is a unique
$a_s \in A$ such that $(a_s \otimes 1) \cdot \Sbas_{00} = s$.

\begin{Theorem}
$\Bimod{S}$ is a weak hopfish antipode with respect to the pairing defined for $\zeta \in \Hom_A(\Bimod{\epsilon}, \Bimod{\Delta})$ and $s \in \Bimod{S}$ by
\begin{equation}
\label{eq:TorusPairing}
 \Pairing{\zeta}{s} :=
 \bigl( \zeta \cdot (a_s \otimes 1) \bigr)^0_0\,,
\end{equation}
where the right hand side denotes the matrix element defined in Eq.~\eqref{eq:MatrixEl}.
\end{Theorem}

\begin{proof}
First of all, we compute, as in the derivation of
Eq.~\eqref{eq:HomAAact}, the action of $A \otimes A$ on the matrix
elements of $\zeta \in \Hom_A(\Bimod{\epsilon}, \Bimod{\Delta})$:
\begin{equation}
\label{eq:HomAAact2}
 \bigl( \zeta \cdot (e_{m_1 j_1} \otimes e_{m_2 j_2})
  \bigr)^l_n
 = \E^{-\I \lambda j_1(n + m_1 - m_2)} \E^{\I m_2 l}
   \zeta^{l+ j_1 + j_2}_{n+m_1-m_2}
 \,.
\end{equation}

We need to show $A \otimes A$-tensoriality of the pairing. Note, that the map $s \mapsto a_s$ is by construction left $A$-linear, $a_{(b \otimes 1) \cdot s} = b * a_s$. Hence,
\[
 \Pairing{\zeta}{(b\otimes 1) \cdot s}
 = \bigl( \zeta \cdot ([b \ast a_s] \otimes 1) \bigr)^0_0
 = \bigl( [\zeta \cdot (b \otimes 1)] \cdot (a_s \otimes 1) \bigr)^0_0
 = \Pairing{\zeta \cdot (b \otimes 1)}{s} \,,
\]
so the pairing is $A \otimes \bbC$-tensorial. Using Eq.~\eqref{eq:HomAAact2} we obtain
\[
 \bigl( \zeta \cdot (S(\Abas_{m,j}) \otimes 1) \bigr)^0_0
 = \bigl( \zeta \cdot (\E^{- \I \lambda m j}
 \Abas_{-m,j} \otimes 1) \bigr)^0_0 \,
 = \zeta^j_{-m}
 = \bigl( \zeta \cdot (1 \otimes \Abas_{mj}) \bigr)^0_0 \,,
\]
which implies
\[
 \Pairing{\zeta \cdot (S(b) \otimes 1)}{\Sbas_{00}}
 = \Pairing{\zeta \cdot (1 \otimes b)}{\Sbas_{00}} \,.
\]
Furthermore, since $\Bimod{S}$ is the modulation of $S$, we have $(1 \otimes b) \cdot s = ([a_s * S(b)] \otimes 1) \cdot \Sbas_{00}$. We conclude
\[
\begin{split}
 \Pairing{\zeta}{(1 \otimes b) \cdot s}
 &= \Pairing{\zeta}{([a_s * S(b)] \otimes 1) \cdot \Sbas_{00}}
 = \Pairing{\zeta \cdot ([a_s * S(b)] \otimes 1)}{\Sbas_{00}} \\
 &= \Pairing{\zeta \cdot (a_s \otimes b)}{\Sbas_{00}}
 = \Pairing{\zeta \cdot (1 \otimes b)}{s}
 \,,
\end{split}
\]
that is, the pairing is $\bbC \otimes A$ tensorial.

It remains to show that the pairing is non-degenerate. Given $\zeta \in \Hom_A(\Bimod{\epsilon}, \Bimod{\Delta})$, assume that $\Pairing{\zeta}{s} =0$ for all $s \in \Bimod{S}$. Then $(\zeta \cdot (a \otimes 1))^0_0 = 0$ for all $a \in A$.  By Eq.~\eqref{eq:HomAAact2} we get
\[
 \bigl( \zeta \cdot (\Abas_{nl} \otimes 1) \bigr)^0_0
 = \E^{-\I \lambda n l} \zeta^{l}_{n} = 0 \,,
\]
for all $n,l \in \mathbb{Z}$. Hence, $\zeta = 0$. Given $s \in \Bimod{S}$ assume now that $\Pairing{\zeta}{s} = 0$ for all $\zeta \in \Hom_A(\Bimod{\epsilon}, \Bimod{\Delta})$. Write $a_s = \sum_{nl} \alpha^{nl} e_{nl}$ with coefficients $\alpha^{nl} \in \mathbb{C}$. Let now $\zeta \in \Hom_A(\Bimod{\epsilon}, \Bimod{\Delta})$ be the map with only one non-zero matrix element $\zeta^{l}_{n} = 1$. Then
\[
 \Pairing{\zeta}{s} = \E^{-\I \lambda n l} \alpha^{nl} = 0 \,.
\]
By choosing such a $\zeta$ for all $n,l \in \mathbb{Z}$ we
conclude that all $\alpha^{nl}$s vanish, so $a = 0$.
\end{proof}

\section{The tensor product of modules}

The coproduct bimodule $\Bimod{\Delta}$ determines a tensor product
operation
$\Motimes$ on right $A$-modules $T$ and $T'$ defined by
\begin{equation}
\label{eq:TensModDef}
  T \Motimes T'
  := (T \otimes T') \otimes_{A \otimes A}
  \Bimod{\Delta} \,.
\end{equation}
By the axioms for a sesquiunital sesquialgebra, the tensor product
 descends to a monoid structure on the isomorphism classes of right
 $A$-modules.  For the tensor product on isomorphism classes, Axiom
 (H1) implies associativity up to isomorphism and Axiom (H2) implies
 that the counit $\Bimod{\epsilon}$ is the unit element.

For simplicity of notation, in this section we shall write the
convolution product of $a,b \in A$ as $ab$, without the star.

\subsection{Simple modules generated by an eigenvector of a unitary}

We now consider a class of simple right $*$-modules of $A$
generated by an eigenvector of a unitary element $u \in U(A)$.
Since the eigenvalue $z$ must be in $U(1)$ we can rescale the
eigenvalue to $1$ by choosing the unitary $z^{-1} u$ instead. The most obvious examples of such modules are obtained as quotient of $A$ by the right ideal $(u-1)A$, which is the smallest possible annihilator of an eigenvector of $u$ with eigenvalue $1$,
\begin{equation}
\label{eq:QuotientModule}
   T = A/ (u-1)A \,.
\end{equation}
We will prove that such a module is simple if and only if $u$ does
not have roots.

\begin{Lemma}
The invertible (unitary) elements of $A$ are the invertible
(unitary) scalar multiples of the basis elements,
\[
    A^\times = \{ \mu\Abas_{pq} \,|\,
                  \mu \in \bbC^\times,\, p,q \in \bbZ \} \,,\quad
    U(A) = \{ \mu\Abas_{pq} \,|\, \mu \in \bbC^\times,\, |\mu|=1,\ p,q \in \bbZ \} \,.
\]
Furthermore, $A$ is a division ring.
\end{Lemma}
\begin{proof}
$A$ is graded with respect to the $\bbZ^2$-grading given by
$\deg(\Abas_{pq}) = (p,q)$. Equip $\bbZ^2$ with the
lexicographical ordering.  For a general nonzero $a \in A$, the maximal and
minimal degrees of its nonzero homogeneous terms are denoted by
$\deg_{\max}(a)$ and $\deg_{\min}(a)$ respectively; the degree of
$0$ is denoted by $\varnothing$ and considered to be less than any
other degree.
Now $\deg_{\max}(a b) = \deg_{\max}(a)
+ \deg_{\max}(b)$ and $\deg_{\min}(a b) = \deg_{\min}(a) +
\deg_{\min}(b)$.  Now let $a$ be invertible with inverse $b$.
Since $\deg(1) = (0,0)$,  $a b = 1$
implies $\deg_{\max}(a) = -\deg_{\max}(b)$ and $\deg_{\min}(a) =
-\deg_{\min}(b)$. Since furthermore $\deg_{\max}(a) \geq
\deg_{\min}(a)$ it follows that $\deg_{\max}(a) = \deg_{\min}(a)$.
Hence, $a$ is of homogeneous degree, that is, proportional to a
basis element, and all basis elements are unitary. By a similar
reasoning it follows from $a b = 0$ and $\deg(0) = \varnothing$
that $a$ and $b$ cannot both be nonzero.
\end{proof}

\begin{Proposition}
\label{th:SimpleMod}
  If $u \in U(A)$ does not have any roots then the right $A$-module $A/(u-1)A$ is simple.
\end{Proposition}

\begin{proof}
To prove that $A/(u-1)A$ is simple, we prove that if $v\notin (u-1)A$,
then $I:=vA+(u-1)A$ is equal to $A$. We start with several
observations.
\begin{enumerate}

\item $u$ can not be a constant because otherwise $u$ has an $n$'th
root for every $n$. Therefore $u=\mu a_{pq}$ with $\mu \in U(1)$, $(p,q)\ne (0,0)$.

\item If $\gcd(p,q) = d>1$, then $\mu a_{pq}=(\sqrt[d]{\mu}\, e^{-\frac{i\lambda
pq(d-1)}{2d^2}}a_{\frac{p}{d}, \frac{q}{d}})^d$. Therefore, $p$
and $q$ have to be relatively prime.
\item Since $\gcd(p,q)=1$, there is $(r,s)\in \bbZ\times \bbZ$ such that
$(p,q)$ and $(r,s)$ form a basis of $\bbZ\times \bbZ$. Therefore,
$a_{pq}$ and $a_{rs}$ generate $A$, and
$v=\sum_{k,l}\alpha_{k,l}a_{pq}^k a_{rs}^l$. Because
$a_{pq}^k-1=(a_{pq}-1)(a_{pq}^{k-1}+\cdots+1)\in (u-1)A$, we can
assume that $v=\sum_l\alpha_{l}a_{rs}^l\ne 0$, a Laurent
polynomial in $a_{rs}$. Furthermore, some $\alpha_l$ is non-zero and $a_{rs}^l$ is
invertible, so $v (a_{rs}^l)^{-1}$ generates the same ideal as $v$ and has a non-zero constant term. Hence, we can assume without loss of generality that $v$ is a polynomial with a non-zero constant term.
\end{enumerate}

Now we show that $I$ is equal to $A$ by proving that $1\in I$.
Assume that $v$ be a polynomial of degree $n$. Since
$a_{pq}a_{rs}^na_{pq}^{-1}=e^{-i\lambda nps}a_{rs}^n$, when $n>0$,
$v'=e^{-i\lambda nps}uvu^{-1}-v=(u-1)e^{-i\lambda
nps}vu^{-1}+v(e^{-i\lambda nps}u^{-1}-1)\in I$ is a polynomial of
degree $n-1$ with a constant term equal to 1. By repeating this
construction, we conclude that 1 is in $I$.
\end{proof}
\begin{Proposition}
\label{th:TDecomp}
If $u \in U(A)$ has a primitive $d$-th root $\sqrt[d]{u}$ then $A/(u-1)A$ can be decomposed as \begin{equation}
\label{eq:TDecomp}
  A/(u-1)A \cong \bigoplus_{n=0}^{d-1} A/(\E^{-\frac{2\pi\I n}{d}}\sqrt[d]{u} -1)A
\end{equation}
into a direct sum of simple right $A$-modules.
\end{Proposition}
\begin{proof}
All primitive $d$-th roots of $u$ can be obtained by multiplying $\sqrt[d]{u}$ with a $d$-th root of unity in $\bbC$, so for every integer $n$, $0\leq n \leq d-1$ we can factorize
\begin{equation}
\label{eq:Decomp1}
  u - 1 =
  \left( \sum_{k=0}^{d-1}
  \bigl( \E^{- \frac{2 \pi \I n}{d}} \sqrt[d]{u} \bigr)^k \right)
  \left( \E^{- \frac{2 \pi \I n}{d}} \sqrt[d]{u} -1 \right) \,.
\end{equation}
Let $\xi := 1 + (u-1)A$ be the canonical cyclic vector of
$A/(u-1)A$. Define for all integers $0\leq n \leq d-1$ the vectors
\[
  \xi_n := \xi \cdot \frac{1}{\sqrt{d}} \sum_{k=0}^{d-1}
  \bigl( \E^{- \frac{2 \pi \I n}{d}} \sqrt[d]{u} \bigr)^k \,.
\]

Consider the cyclic submodule $\xi_n \cdot A$ which is isomorphic
to $A/ \Ann(\xi_n)$. By Eq.~\eqref{eq:Decomp1} the annihilator of
$\xi_n$ contains the ideal $I := \bigl( \E^{- \frac{2 \pi \I
n}{d}} \sqrt[d]{u} -1 \bigr)A$. Since $\sqrt[d]{u}$ is a primitive
root, it does not have proper roots itself, so it follows from
Proposition~\ref{th:SimpleMod} that $I$ is maximal and, hence,
equal to the annihilator of $\xi_n$. We conclude that
\[
  \xi_n \cdot A \cong A/(\E^{- \frac{2 \pi \I n}{d}} \sqrt[d]{u} - 1)A \,.
\]
The cyclic vector can be retrieved as $\xi := \frac{1}{\sqrt{d}} \sum_{n=0}^{d-1} \xi_n$, which implies that the $\xi_n$'s generate the whole module,
\[
  A/(u-1)A = \xi_0 \cdot A + \ldots + \xi_{d-1} \cdot A \,.
\]
It remains to show that this sum is direct.
Since $\{\Abas_{jk} \,|\, j,k \in \bbZ\}$ is a basis of $A$, the set
of vectors $\{\xi_n \cdot \Abas_{jk} \,|\, j,k \in \bbZ\}$ spans
$\xi_n \cdot A$. The primitive root is of the form $\sqrt[d]{u} =
\E^{-\I\alpha} \Abas_{pq}$. Using Eq.~\eqref{eq:AProduct} for the
multiplication in $A$ we get $(\xi_n \cdot \Abas_{jk}) \cdot
\sqrt[d]{u} = \E^{\I[ \frac{2 \pi n}{d} + \lambda(jq - kp)]} (\xi_n
\cdot \Abas_{jk})$. Hence, each summand $\xi_n \cdot A$ can be
decomposed into eigenspaces of $\sqrt[d]{u}$ with eigenvalues of the
form $\E^{\I[ \frac{2 \pi n}{d} + \lambda k]}$ for some integer $k$.
Because $\frac{\lambda}{2\pi}$ is irrational, the eigenvalues of the
eigenspace decomposition of one summand $\xi_n \cdot A$ are all
different from those of every other summand $\xi_{n'} \cdot A$, $n'
\neq n$. We conclude that the subspaces $\xi_0 \cdot A, \ldots,
\xi_{d-1} \cdot A$ are linearly independent, so their sum is direct.
\end{proof}

\begin{Corollary}
The right $A$-module $A/(u-1)A$ is simple if and only if $u =
\E^{-\I\alpha} \Abas_{pq}$ for $\alpha \in \bbR$ and $p,q \in
\bbZ$ relatively prime.
\end{Corollary}
\begin{proof}
From Eq. \eqref{eq:AProduct} we derive for all $d, p', q' \in
\bbZ$ and $p = d p'$, $q = d q'$:
\begin{equation}
\label{eq:AbasPowers}
  (\Abas_{p'q'})^d = \E^{\frac{\I \lambda p'q' d(d-1)}{2} } \Abas_{dp',dq'}
  \quad\Leftrightarrow\quad
  \sqrt[d]{\Abas_{pq}} = \E^{-\frac{\I \lambda pq(d-1)}{2d^2} }
  \Abas_{\frac{p}{d},\frac{q}{d}} \,.
\end{equation}

This shows that $\Abas_{pq}$ does not have roots if and only if
$p$ and $q$ are relatively prime.
\end{proof}

\begin{Definition}
Let $\alpha \in \bbR$, $p,q \in \bbZ$ relatively prime. Define
\[
  T^{\alpha}_{pq} := A/(\E^{-\I\alpha} \Abas_{pq} - 1)A \,.
\]
\end{Definition}

\subsection{Construction of a basis}

We will find a basis for the simple modules
$T^\alpha_{pq}$. Let $u := \E^{-\I\alpha}\Abas_{pq}$ and let $\xi
:= 1 + (u-1)A$ be the canonical cyclic vector. The vectors
$\xi'_{jk} := \xi \cdot \Abas_{jk}$ span $T^\alpha_{pq}$. Using
Eq.~\eqref{eq:AProduct}, we get on the one hand,
\begin{equation}
\label{eq:uEigen}
  \xi'_{jk} \cdot u = \E^{\I\lambda(j q - k p)} \xi'_{jk} \, ,
\end{equation}
on the other hand, $\xi'_{jk} \cdot u
  = \xi \cdot \bigl(
    \E^{\I(-\alpha + \lambda jq)} \Abas_{j+p,k+q} \bigr)
  = \E^{\I(-\alpha + \lambda jq)} \xi_{j+p,k+q}$, which implies
\[
  \xi'_{j+p,k+q} = \E^{\I(\alpha - \lambda k p)} \xi'_{jk} \,.
\]
We need to consider two cases separately.

\textbf{Case 1:} $p \neq 0$. We can rescale the vectors $\xi'_{jk}$ to
\[
  \xi_{jk} :=
  \E^{ \frac{\I}{p}\left\{ - \alpha j + \lambda [ j k p - \frac{1}{2} j(j+p)q ]
  \right\} } \xi'_{jk} \,,
\]
which satisfy $\xi_{jk} = \xi_{j+p,k+q}$. Hence, the
vector $\xi_{jk}$ can be labeled by $[j,k] := (j,k) + (p,q)\bbZ$.
Since $p$ and $q$ are relatively prime, there is a bijection
\begin{equation}
\label{eq:LabelDef}
  \nu: \bbZ^2/(p,q)\bbZ \stackrel{\cong}{\longrightarrow} \bbZ \,,\qquad
  [j,k] \longmapsto jq - kp \,.
\end{equation}
Eq.~\eqref{eq:uEigen} implies that, if $[j,k] \neq [j',k']$ then $\xi_{[j,k]}$ and $\xi_{[j',k']}$ are eigenvectors of $u$ with different eigenvalues. We conclude that the set of vectors
\begin{equation}
\label{eq:Tbasis}
  B := \{ \xi_{[j,k]} \,|\, [j,k] \in \bbZ^2/(p,q)\bbZ \}
\end{equation}
is an orthonormal basis of $T^{\alpha}_{pq}$. The action of the generators of $A$ on this basis is
\begin{equation}
\label{eq:Tact1a}
\begin{aligned}
  \xi_{[j,k]} \cdot \Abas_{10}
  &= \E^{ \frac{\I}{p}\left[ \alpha
     + \lambda ( j q - k p) + \frac{1}{2} \lambda q(p+1) \right]} \xi_{[j+1,k]} \,,\\
  \xi_{[j,k]} \cdot \Abas_{01}
  &= \xi_{[j,k+1]} \,.
\end{aligned}
\end{equation}
If we label the basis instead by $\bbZ$ via the bijection defined in Eq.~\eqref{eq:LabelDef}, $\xi_n := \xi_{\nu^{-1}(n)}$, the action reads
\begin{equation}
\label{eq:Tact1b}
\begin{aligned}
  \xi_{n} \cdot \Abas_{10}
  &= \E^{ \frac{\I}{p}\left[ \alpha
     + \lambda n + \frac{1}{2} \lambda q(p+1) \right]}
     \xi_{n+q} \,,\\
  \xi_{n} \cdot \Abas_{01}
  &= \xi_{n-p} \,.
\end{aligned}
\end{equation}

\textbf{Case 2:} $p = 0$. Because the module $T^\alpha_{0q}$ is
assumed to be simple, $q$ must be nonzero. We can rescale
\[
  \xi_{jk} := \E^{ -\frac{\I \alpha k}{q} } \xi_{jk} \,,
\]
such that $\xi_{jk} = \xi_{j+p,k+q}$. The set of vectors defined
in Eq.~\eqref{eq:Tbasis} is again an orthonormal basis. The action
of the generators of $A$ on this basis is
\begin{equation}
\label{eq:Tact2a}
\begin{aligned}
  \xi_{[j,k]} \cdot \Abas_{10}
  &=      \xi_{[j+1,k]} \,,\\
  \xi_{[j,k]} \cdot \Abas_{01}
  &= \E^{ \frac{\I}{q} \left[ \alpha + \lambda j q \right]}
  \xi_{[j,k+1]} \,.
\end{aligned}
\end{equation}
Again, labeling the basis by $\bbZ$, we read the action to be
\begin{equation}
\label{eq:Tact2b}
\begin{aligned}
  \xi_{n} \cdot \Abas_{10}
  &=      \xi_{n+q} \,,\\
  \xi_{n} \cdot \Abas_{01}
  &= \E^{ \frac{\I}{q} \left[ \alpha + \lambda n \right]}
  \xi_{n} \,.
\end{aligned}
\end{equation}

We can use these formulas to construct modules in the case that $p$ and $q$ are not relatively prime:

\begin{Definition}
\label{th:Tdef}
Let $T$ be the inner product space spanned by the orthonormal basis $\{\xi_n \,|\, n \in \bbZ \}$. Let $\alpha \in \bbR$, $p,q \in \bbZ$ where $p \neq 0$ or $q \neq 0$. The right action of $A$ on $T$ given by Eqs.~\eqref{eq:Tact1b} for $p \neq 0$ and by Eqs.~\eqref{eq:Tact2b} for $p = 0$ defines a right $A$ $*$-module, which we denote by $T^\alpha_{pq}$.
\end{Definition}

\subsection{Isomorphism classes of the modules}

We conclude the general study of the modules $T^\alpha_{pq}$ by giving
a criterion for two modules to be isomorphic.

\begin{Proposition}
The modules $T^{\alpha}_{pq}$ and $T^{\beta}_{rs}$ are isomorphic
iff $(p,q) = \pm (r,s)$ and $\alpha = \pm \beta + n \lambda$ for
some $n \in \bbZ$.
\end{Proposition}
\begin{proof}
First, let us assume that $T^{\alpha}_{pq}$ and $T^{\beta}_{rs}$ are simple.
By Eq. \eqref{eq:Tact1b}, the matrix with respect to the basis
$\{\xi_n\}$ for $a_{rs}$ acting on
$T^{\alpha}_{pq}$ has all its nonzero elements on a diagonal which is
shifted from the main diagonal by $rq-sp$ units.  Unless $rq-sp=0$,
this matrix has no eigenvectors.  Since all the $\xi_n$ are
eigenvalues for the action of $a_{pq}$ on $T^{\alpha}_{pq}$, it
follows that the two modules in question can be isomorphic if and only
if the integer vectors $(p,q)$ and $(r,s)$ are collinear.  Now
$\gcd(p,q)=\gcd(r,s)=1$ implies that $(p,q)=\pm(r,s)$.  The statement
about $\alpha$ and $\beta$ follows from a comparison of the
eigenvalues.  This completes the proof of the ``only if'' part of our
proposition.  The ``if'' part follows immediately from the definition
of the modules.

Now if $\Abas_{pq}$ has a primitive $d$-th root,  $T^{\alpha}_{pq}$ is the direct sum of the modules generated by each of $\xi_k,\xi_{k+1}, \ldots, \xi_{k+d-1}$ for some $k \in \bbZ$. For convenience we choose $k := - \frac{pq(d-1)}{2d}$. Labelling the basis of $\xi_{k + l} \cdot A$ as $\eta_j := \xi_{k + l + d j}$, we can read off Eqs.~\eqref{eq:Tact1b} and \eqref{eq:Tact2b} that $\xi_{k+l} \cdot A \cong T^{[(\alpha + \lambda l)/d]}_{p/d,q/d}$. We thus obtain the decomposition into simple modules
\begin{equation*}
  T^{\alpha}_{pq} \cong \bigoplus_{l=0}^{d-1}
  T^{[\frac{\alpha + \lambda l}{d}]}_{\frac{p}{d}, \frac{q}{d}} \,.
\end{equation*}
For $T^{\alpha}_{pq}$ and $T^{\beta}_{rs}$ to be isomorphic the simple modules of this decomposition have to be pairwise isomorphic. This is the case if there are $0 \leq l,m < d$ and $n \in \bbZ$ such that $(\alpha + \lambda l)/d = \pm (\beta + \lambda m)/d + \lambda n \Leftrightarrow \alpha = \pm \beta + \lambda n'$, where $n' = d n - l \pm m$.
\end{proof}

\subsection{Calculation of the tensor product}

We proceed to calculate the tensor product of two modules
$T^{\alpha_1}_{p_1 q_1}$ and $T^{\alpha_1}_{p_1 q_1}$ of
Definition~\ref{th:Tdef}. As a vector space we have
\[
 T^{\alpha_1}_{p_1 q_1} \Motimes T^{\alpha_2}_{p_2 q_2}
 = \Span \{ \xi^{1}_{k_1} \otimes \xi^{2}_{k_2}
 \otimes \Dbas_{n_1 n_2 l} \,|\,
 k_1,k_2,n_1,n_2,l \in \mathbb{Z} \} / V \,,
\]
where $\xi_{k_1}^{1}$ and $\xi_{k_2}^{2}$ are the basis vectors of
$T^{\alpha_1}_{p_1 q_1}$ and $T^{\alpha_1}_{p_1 q_1}$ from
Definition~\ref{th:Tdef}, and where $V$ is the ideal (the vector
space) generated by the relations
\[
 (\xi^{1}_{k_1} \cdot a_1) \otimes
 (\xi^{2}_{k_2} \cdot a_2) \otimes
 \Dbas_{n_1 n_2 l}
 =
 \xi^{1}_{k_1} \otimes \xi^{2}_{k_2} \otimes
 ( a_1 \otimes a_2) \cdot \Dbas_{n_1 n_2 l} \,,
\]
for all $a_1, a_2 \in A$.  Using $\Dbas_{n_1 n_2 l} = \E^{-\I
\lambda n_1 l} (\Abas_{n_1,-l} \otimes \Abas_{n_2 0}) \cdot
\Dbas_{000}$ we obtain that $T^{\alpha_1}_{p_1 q_1} \Motimes
T^{\alpha_2}_{p_2 q_2}$ is spanned by vectors of the form
\[
 \xi'_{k_1, k_2} :=
 \xi^{1}_{k_1} \otimes \xi^{2}_{k_2} \otimes \Dbas_{0 0 0} \,.
\]
These vectors are not yet linearly independent.  We still have to
consider the relation generated by the action of $ 1 \otimes
\Abas_{0 1}$. Using $(1 \otimes \Abas_{01}) \cdot \Dbas_{000} =
\Dbas_{00,-1} = (\Abas_{01} \otimes 1) \cdot \Dbas_{000}$, we
obtain the relation
\begin{equation}
\label{eq:TensorRel0}
 \xi^{1}_{k_1} \otimes (\xi^{2}_{k_2} \cdot \Abas_{0 1})
 \otimes \Dbas_{0 0 0}
 = (\xi^{1}_{k_1} \cdot \Abas_{0 1}) \otimes \xi^{2}_{k_2}
 \otimes \Dbas_{0 0 0} \,.
\end{equation}
We have to distinguish three cases:

\textbf{Case 1:} $p_1 \neq 0$ and $p_2 \neq 0$. Under this assumption relation~\eqref{eq:TensorRel0} becomes
\[
 \xi'_{k_1 + p_1 , k_2 - p_2} = \xi'_{k_1, k_2} \,.
\]
We conclude that we can label these vectors uniquely by $\bbZ^2/(p_1,-p_2)\bbZ$,
\[
  \xi_{[k_1,k_2]} := \xi'_{k_1,k_2} \,,\qquad
  [k_1,k_2] := (k_1,k_2)+ (p_1,-p_2)\bbZ \,.
\]
By construction, these vectors form a basis of the tensor product
module. The action of $\Abas_{10} \in A$ on this basis is given by
\[
\begin{split}
  \xi_{[k_1, k_2]} \cdot \Abas_{10}
  &= \xi^{1}_{k_1} \otimes \xi^{2}_{k_2}
  \otimes [ \Dbas_{0 0 0} \cdot \Abas_{10} ] \\
  &= \xi^{1}_{k_1} \otimes \xi^{2}_{k_2}
  \otimes [(\Abas_{10} \otimes \Abas_{10}) \cdot \Dbas_{0 0 0}  ] \\
  &= ( \xi^{1}_{k_1} \cdot \Abas_{10} ) \otimes
     ( \xi^{2}_{k_2} \cdot \Abas_{10} )
     \otimes \Dbas_{0 0 0} \,,
\end{split}
\]
and similarly for $\Abas_{01}$. Inserting~\eqref{eq:Tact1b} yields
\begin{equation}
\label{eq:TotimesTact1}
\begin{aligned}
  \xi_{[k_1, k_2]} \cdot \Abas_{10}
  &= \E^{\I\left\{
     \frac{\alpha_1 + \lambda[k_1 + q_1(p_1+1)/2]}{p_1}
    +\frac{\alpha_2 + \lambda[k_2 + q_2(p_2+1)/2]}{p_2} \right\}}
  \xi_{[k_1 + q_1, k_2 + q_2]} \\
  \xi_{[k_1, k_2]} \cdot \Abas_{01}
  &= \xi_{[k_1 - p_1, k_2]} \,.
\end{aligned}
\end{equation}
Using the bijection
\[
\begin{aligned}
 \nu: \mathbb{Z}^2/(p_1,-p_2)\bbZ &\stackrel{\cong}{\longrightarrow}
 (\mathbb{Z}/ \gcd(p_1,p_2) \mathbb{Z}) \times \bbZ\\
 [k_1,k_2] &\longmapsto
 \Bigl( q_2 k_1 - q_1 k_2 \bmod \gcd(p_1,p_2)\bbZ ,
       \frac{p_2 k_1 + p_1 k_2}{\gcd(p_1, p_2)} \Bigr) \,,
\end{aligned}
\]
we relabel the basis by setting $\xi^{(m)}_n :=
\xi_{\nu^{-1}(m,n)}$. The action~\eqref{eq:TotimesTact1} then
takes the form
\begin{equation}
\label{eq:TotimesTact2}
\begin{aligned}
  \xi^m_{n} \cdot \Abas_{10}
  &= \E^{ \frac{\I}{p} \left[ \alpha + \lambda n + \frac{1}{2} q(p+1) \right] }
  \xi^m_{n + q} \\
  \xi^m_{n} \cdot \Abas_{01}
  &= \xi^m_{n - p} \,,
\end{aligned}
\end{equation}
where
\[
  p := \lcm(p_1,p_2) \,,~
  q := \frac{p_1 q_2 + p_2 q_1}{\gcd(p_1, p_2)} \,,~
  \alpha := \frac{\alpha_1 p_2 + \alpha_2 p_1}{\gcd(p_1,p_2)}
       + \frac{\lambda p(q_1 + q_2 - q)}{2} \,.
\]
Comparing this action with~\eqref{eq:Tact1b} we infer
\begin{equation}
\label{eq:TotimesT1}
 T^{\alpha_1}_{p_1 q_1} \Motimes T^{\alpha_2}_{p_2 q_2}
 \cong \gcd(p_1,p_2)\, T^{\alpha}_{p q} \,,
\end{equation}
where the prefactor on the right hand side denotes the direct sum of $\gcd(p_1,p_2)$ copies of the same module.

\textbf{Case 2:} either $p_1 = 0$ or $p_2 = 0$. The calculation of the tensor product for this case is very similar to the preceding case. The result is again given by Eq.~\eqref{eq:TotimesT1}.

\textbf{Case 3:} $p_1 = 0$ and $p_2 = 0$. This case is quite different. Using Eq.~\eqref{eq:Tact2b} for the action of $\Abas_{01}$, relations~\eqref{eq:TensorRel0} become
\[
 \xi^{1}_{k_1} \otimes
 (\E^{ \frac{\I}{q_2} \left[ \alpha_2 + \lambda k_2 \right]} \xi^{2}_{k_2})
 \otimes \Dbas_{0 0 0} \\
 = (\E^{ \frac{\I}{q_1} \left[ \alpha_1 + \lambda k_1 \right]} \xi^{1}_{k_1})
   \otimes \xi^{2}_{k_2} \otimes \Dbas_{0 0 0} \,,
\]
which is equivalent to
\[
  \Bigl(
   \E^{ \frac{\I}{q_2} \left[ \alpha_2 + \lambda k_2 \right]}
  -\E^{ \frac{\I}{q_1} \left[ \alpha_1 + \lambda k_1 \right]}
  \Bigr) \xi'_{k_1,k_2} = 0 \,.
\]
It follows that $\xi'_{k_1,k_2}$ is zero unless
\[
  \frac{\alpha_1 q_2 - \alpha_2 q_1 + \lambda (k_1 q_2 - k_2 q_1)}{q_1 q_2}
  \cong 0 \bmod 2\pi\bbZ \,.
\]

Assuming (without loss of generality) that $0 \leq
\alpha_1,\alpha_2 < 2\pi$ and using that $\frac{\lambda}{2\pi}$ is
irrational, we conclude that all $\xi'_{k_1,k_2}$ vanish unless
\[
  r := - \frac{\alpha_1 q_2 - \alpha_2 q_1}{\lambda \gcd(q_1,q_2)}
\]
is an integer modulo integer multiples of $\lcm(q_1,
q_2)\frac{2\pi}{\lambda}$.

In this case a basis for the tensor module is given by
\[
  B := \Bigl\{ \xi'_{k_1,k_2} \,\Big|\, k_1,k_2 \in \bbZ \,,
  \frac{k_1 q_2 - k_2 q_1}{\gcd(q_1,q_2)} = r \Bigr\} \,.
\]
We can relabel the basis by $\bbZ$: Let $s_1, s_2$ be integers such that
\begin{equation}
\label{eq:sChoice}
  \frac{s_1 q_2 - s_2 q_1}{\gcd(q_1,q_2)} = 1 \,.
\end{equation}
This can be used to construct a bijection
\[
\begin{aligned}
  \nu: \bbZ &\stackrel{\cong}{\longrightarrow}
    \Bigl\{(k_1,k_2) \in \bbZ \,\Big|\,
    \frac{k_1 q_2 - k_2 q_1}{\gcd(q_1,q_2)} = r \Bigr\} \,,\\
  n & \longmapsto r ( s_1, s_2) + \frac{k}{\gcd(q_1,q_2)} (q_1,q_2) \,,
\end{aligned}
\]
and relabel the basis by $\xi_n := \xi'_{\nu(n)}$. In terms of the
relabeled basis the action of the generators $\Abas_{10},
\Abas_{01}$ takes the form of Eqs.~\eqref{eq:Tact2b} with
\[
  q := \gcd(q_1,q_2) \,, \quad
  \alpha := s_1 \alpha_2 - s_2 \alpha_1 \,.
\]

For convenience we summarize the results obtained in this section.

\begin{Theorem}
\label{th:TensorProduct} Let $T^{\alpha_1}_{p_1 q_1}$ and
$T^{\alpha_2}_{p_2 q_2}$ be the right $A$-modules of
Definition~\ref{th:Tdef}, their tensor product being defined in
Eq.~\eqref{eq:TensModDef}.

For $p_1 \neq 0$ or $p_2 \neq 0$ we have:
\[
 T^{\alpha_1}_{p_1 q_1} \Motimes T^{\alpha_2}_{p_2 q_2}
  \cong \gcd(p_1,p_2)\, T^{\alpha}_{p q} \,,
\]
where
\begin{equation}
\label{eq:alphapqDef}
  p := \lcm(p_1,p_2) \,,~
  q := \frac{p_1 q_2 + p_2 q_1}{\gcd(p_1, p_2)} \,,~
  \alpha := \frac{\alpha_1 p_2 + \alpha_2 p_1}{\gcd(p_1,p_2)} \,.
\end{equation}
For $p_1 = 0$ and $p_2 = 0$ we have:
\[
  T^{\alpha_1}_{0, q_1} \Motimes T^{\alpha_2}_{0, q_2} \cong
  \begin{cases}
  T^{\alpha}_{0, q} \,, &\mathrm{for}\quad
  \frac{\alpha_1 q_2 - \alpha_2 q_1}{\lambda \gcd(q_1,q_2)} \in \bbZ \bmod \lcm(q_1,q_2)\frac{2\pi}{\lambda}\\
  0 \,,& \mathrm{otherwise}
  \end{cases} \,,
\]
where
\begin{equation}
\label{eq:alphapqDef2}
q := \gcd(q_1,q_2) \,, \quad
  \alpha := s_1 \alpha_2 - s_2 \alpha_1
  \,,\quad s_1,s_2 \in \bbZ \,:\quad
  \frac{s_1 q_2 - s_2 q_1}{\gcd(q_1,q_2)} = 1 \,.
\end{equation}
\end{Theorem}

We end this paper with the following remarks.

\begin{Remark}
We have dropped the term $\frac{\lambda p(q_1 + q_2 - q)}{2}$ in
the expression (\ref{eq:alphapqDef}) for $\alpha$ because
$p(q_1+q_2-q)$ is always even.
\end{Remark}

\begin{Remark}
The formulas for $\Motimes$ in Theorem \ref{th:TensorProduct} have
the following simple interpretation as addition of fractions, e.g.
\[
\frac{p_1}{q_1}+\frac{p_2}{q_2}=\frac{p}{q}.
\]
This is connected with the second symplectic groupoid structure
mentioned in the Introduction.
\end{Remark}

\begin{Remark}
The operation $\Motimes$ extends to an associative commutative
product on the free abelian group $\mathcal {R}$ generated by the
$T_{pq}^\alpha$, with the module $T_{10}^{0}=\bepsilon$ as unit of
$\mathcal {R}$. The resulting unital ring $\mathcal {R}$ seems to
contain the necessary information to reconstruct the group
structure on the original quotient $S^1 / \bbZ$.
\end{Remark}


\begin{thebibliography}{99}

\bibitem{ba-la:2-groups}
Baez, John C.,  and Lauda, Aaron D., Higher-dimensional algebra. {V}. 2-groups, {\em Theory Appl. Categ.} {\bf 12} (2004), 423--491

\bibitem{be:stack}
Behrend, K., Cohomology of stacks.  {\em Intersection theory and moduli},  (electronic), {\em ICTP Lect. Notes}, XIX, Abdus Salam Int. Cent. Theoret. Phys., Trieste, 2004, 249--294

\bibitem{be-xu:stack}
Behrend, K., and Xu, P., Differentiable stacks and gerbes, math.DG/0605694 (2006).

\bibitem{bo-fl-ge-pi:hidden}
Bonneau, P., Flato, M., Gerstenhaber, M., and Pinczon, G., The hidden group structure of quantum groups: strong duality, rigidity and preferred deformations, {\em Comm. Math. Phys.} {\bf 161} (1994), 125--156.

\bibitem{ka-ra:para-hopf}
Khalkhali, M., and Rangipour, B., Para-Hopf algebroids and their
cyclic cohomology, {\em Lett. Math. Phys. }70 (2004), no. 3,
259--272.

\bibitem{lu-we:groupoides}
Lu, J.-H., and Weinstein, A., Groupo\"{\i}des symplectiques doubles
des groupes de Lie-Poisson, {\em C. R. Acad. Sci. Paris} {\bf 309}
(1989), 951--954.

\bibitem{mr:stability}
Mr{\v{c}}un, J. {\em Stablility and invariants of
{H}ilsum-{S}kandalis maps}, PhD thesis, Utrecht University,
Utrecht, 1996, math.DG/0506484.

\bibitem{ri:quantization}
Rieffel, M., {\em Deformation quantization for actions of $R\sp
d$}, Mem. Amer. Math. Soc. 106 (1993), no. 506.

\bibitem{ta-we-zh:hopfish}
Tang, X., Weinstein, A., and Zhu, C., Hopfish algebras, to appear in
{\em Pacific   J. Math}, math.QA/0510421.

\bibitem{tz:stack}
Tseng, H., and Zhu. C., Integrating {L}ie algebroids via stacks.
{\em Compositio Mathematica}, Volume 142 (2006), 251--270.

\bibitem{we:rotation}
Weinstein, A., Symplectic groupoids, geometric quantization,
and irrational rotation algebras, {\em Symplectic geometry, groupoids,
and integrable systems, S\'{e}minaire sud-Rhodanien
de g\'{e}om\'{e}trie \`{a} Berkeley (1989)}, P. Dazord and A.
Weinstein, eds., Springer-MSRI Series (1991), 281--290.

\bibitem{wo:quant-gp}
Woronowicz, S., Twisted ${\rm SU}(2)$ group, An example of a
noncommutative differential calculus, {\em Publ. Res. Inst. Math.
Sci.}, 23 (1987), 117--181.

\bibitem{za:quantum}
Zakrzewski, S., Quantum and classical pseudogroups, I and II,
{\em  Comm. Math. Phys.} {\bf 134} (1990), 347-370, 371-395.

\end{thebibliography}
\end{document}